# Identities for number series and their reciprocals: Dirac delta function approach

S. M. Abrarov [1], R. M. Abrarov [2]

April 20, 2007

**Abstract**

Dirac delta function (delta-distribution) approach can be used as efficient method to derive identities for number series and their reciprocals. Applying this method, a simple proof for identity relating prime counting function ($\pi$-function) and logarithmic integral ($Li$-function) can be obtained.

## I. Introduction

Since Paul Dirac introduced the delta function [1], it remains a versatile mathematical tool in many applications of the modern sciences and information technologies. Being a physicist P. Dirac intuitively believed in the existence of the delta function as a true mathematical object and successfully applied it developing the fundamentals of the relativistic quantum mechanics in his classic works [1, 2].

Nowadays Dirac delta function is widely used in the pure and applied mathematics and engineering physics [3-8]. Particularly the sampling methods, based on Dirac delta function, proved their reliability and are commonly implemented in the systems of telecommunications and signal processing [8].

In this work we introduce some useful applications of delta-distribution in the theory of numbers. Compared to Stieltjes integration methods [9], Dirac delta function approach is more natural and obvious for understanding and can be successfully used in the derivations of the identities for the number series and their reciprocals.

## II. Some properties of $\delta$-function

We will apply two main properties of the $\delta$-function. The first one is the sampling property

$$f(x_0) = \int_{-\infty}^{\infty} f(x)\delta(x-x_0)dx, \qquad (1)$$

while the second one relates itself with Heaviside step function [10]

$$\chi(x-x_0) = \begin{cases} 0, & x < x_0 \\ 1, & x \geq x_0 \end{cases}$$

via integration

$$\chi(x-x_0) = \int_{-\infty}^{\infty} \delta(x-x_0)dx. \qquad (2)$$

It should be noted that integral of the delta function can be determined by infinitely narrow area in the vicinity of $x_0$. Therefore (1) and (2) can be rewritten as

$$f(x_0) = \int_{x_0-\varepsilon}^{x_0} f(x)\delta(x-x_0)dx$$

and

$$\chi(x-x_0) = \int_{x_0-\varepsilon}^{x_0} \delta(x-x_0)dx,$$

respectively, where $\varepsilon$ is infinitely small positive value.

## III. Step function series

### *3.1 Number of pulses*

Consider a set consisting of $N$ numbers $q_1$, $q_2$, $q_3$, ... , $q_N$. Define the series as a sum of reciprocals $q_i$ as follows

$$h(x) = \begin{cases} \sum_{q_i \leq x} \dfrac{1}{q_i}, & x \geq q_{min} \\ 0, & x < q_{min} \end{cases} \qquad (3)$$

where $q_{min}$ is the smallest number in this set. According to the sampling property (1) of the delta-distribution, the series (3) can be rewritten as

$$h(x) = \int_{q_{\min}-\varepsilon}^{x} \frac{1}{y} \sum_{q_i \leq x} \delta_i(y - q_i) \, dy,$$

which, in turn, can be expressed in the form

$$xh(x)' = \sum_{q_i \leq x} \delta_i(x - q_i). \tag{4}$$

In compliance with (2), integrating of (4) yields

$$\int_{q_{\min}-\varepsilon}^{x} yh(y)' \, dy = \sum_{q_i \leq x} \chi_i(x - q_i). \tag{5}$$

RHS of (5) counts the number of stepwise pulses, therefore:

$$\int_{q_{\min}-\varepsilon}^{x} yh(y)' \, dy = \underbrace{1+1+1+\ldots}_{N \text{ times}} = \lfloor x \rfloor,$$

where symbol $\lfloor \, \rfloor$ denotes the floor operator. Equation (5) contains derivative, which may be inconvenient for computing. In order to avoid it, we apply integration by part. This leads to

$$N = xh(x) - \int_{q_{\min}}^{x} h(y) \, dy, \tag{6}$$

where $N$ is the number of the stepwise pulses within interval from $q_{\min}$ up to $x$.

## 3.2 Sum of numbers

Multiplying $x$ to both parts of (4) and integrating the result, we can find that

$$\int_{q_{\min}-\varepsilon}^{x} y^2 h(y)' = \sum_{q_i \leq x} q_i \chi_i(x - q_i). \tag{7}$$

RHS of (7) represents the sum of all numbers up to $x$. Therefore this equation we can be rewritten as

$$\sum_{q_i \leq x} q_i = \int_{q_{\min}-\varepsilon}^{x} y^2 h(y)' \, dy. \tag{8}$$

The derivative in (8) can be excluded again. This provides the following relation

$$\sum_{q_i \leq x} q_i = x^2 h(x) - 2 \int_{q_{\min}}^{x} yh(y) \, dy.$$

## 3.3 Power series

More generally, (4) can be extended for any power of $x$ as given by

$$x^{k+1}h(x)' = x^k \sum_{q_i \leq x} \delta_i(x - q_i),$$

where $k$ is any number. This leads to the relation of the power series

$$\sum_{q_i \leq x} q_i^k = \int_{q_{\min} - \varepsilon}^{x} y^{k+1} h(y)' \, dy. \tag{9}$$

Excluding the derivative in (9), we have

$$\sum_{q_i \leq x} q_i^k = x^{k+1} h(x) - (k+1) \int_{q_{\min}}^{x} y^k h(y) \, dy. \tag{10}$$

## 3.4 Sum of reciprocals

Define the following function as

$$g(x) = \begin{cases} \sum_{q_i \leq x} q_i, & x \geq q_{\min} \\ 0, & x < q_{\min} \end{cases}$$

By analogy with (4), the superposition of the delta functions is given by

$$\frac{g(x)'}{x} = \sum_{q_i \leq x} \delta_i(x - q_i). \tag{11}$$

Dividing both parts of (11) by $x^k$ and integrating the result yields

$$\int_{q_{\min} - \varepsilon}^{x} \frac{g(y)'}{y^{k+1}} \, dy = \sum_{q_i \leq x} \frac{1}{q_i^k} \chi_i(x - q_i) = \sum_{q_i \leq x} \left(\frac{1}{q_i}\right)^k \tag{12}$$

Finally, integrating by part (12) results to

$$\sum_{q_i \leq x} \left(\frac{1}{q_i}\right)^k = \frac{g(x)}{x^{k+1}} + (k+1) \int_{q_{\min}}^{x} \frac{g(y)}{y^{k+2}} \, dy. \tag{13}$$

## IV. Natural and prime numbers

### *4.1 Natural numbers*

We will use the relations derived above for natural and prime numbers. Consider the harmonic number (further index $i$ will be omitted)

$$H_n(x) = \begin{cases} \sum_{n=1}^{\lfloor x \rfloor} \frac{1}{n}, & x \geq 1 \\ 0, & x < 1 \end{cases} \quad (14)$$

where $n = 1, 2, 3 \ldots$ . Substituting (14) into (10), we can find that

$$\sum_{n=1}^{\lfloor x \rfloor} n^k = 1^k + 2^k + 3^k + 4^k + \ldots = x^{k+1} H_n(x) - (k+1) \int_1^x y^k H_n(y) dy.$$

There are two interesting cases. For $k = 0$, we obtain an equation counting the number of pulses, similar to that of (6). However, as the gap between two closest integers is always unity, the number of pulses must be equal to the largest integer $n_{max} \in n$, i.e.:

$$n_{max} = x H_n(x) - \int_1^x H_n(y) dy.$$

For $k = 1$, we have the sum of the arithmetic progression of the natural numbers

$$\sum_{i=1}^{\lfloor x \rfloor} n = 1 + 2 + 3 + 4 + \ldots = x^2 H_n(x) - 2 \int_1^x y H_n(y) dy.$$

Therefore using the formula $(n_{max} + n_{max}^2)/2 = 1 + 2 + 3 + \ldots + n_{max}$, we can write

$$\frac{n_{max} + n_{max}^2}{2} = x^2 H_n(x) - 2 \int_1^x y H_n(y) dy.$$

To derive a relation for harmonic number, we use (13) and substitute $k = 1$. This provides the following identity

$$H_n(x) = \frac{\lfloor x \rfloor + \lfloor x \rfloor^2}{2x^2} + \int_1^x \frac{\lfloor y \rfloor + \lfloor y \rfloor^2}{y^3} dy. \quad (15)$$

It should be noted that the algorithm, built on the basis of (15), significantly accelerates the computation of the harmonic number, especially for the large values of $x$. This is possible to achieve since RHS of (15) includes the set of integers instead of reciprocals. The computation performed directly through (14) involves rapidly decreasing terms as $x$ increases. As a result of the multiple process of division, the calculation becomes problematic due to the decrease of the significant digits after floating point. The use of integers effectively resolves this problem.

## *4.2 Prime numbers*

Define a sum of the reciprocal primes in the following form

$$H_p(x) = \begin{cases} \sum_{p \leq x} \dfrac{1}{p}, & x \geq 2 \\ 0, & x < 2 \end{cases}, \qquad (16)$$

where $p = 2, 3, 5, 7 \ldots$, i.e. all primes smaller than $x$. Substituting (16) into (10), we get

$$\sum_{p \leq 2} p^k = 2^k + 3^k + 5^k + 7^k + \ldots = x^{k+1} H_p(x) - (k+1) \int_2^x y^k H_p(y) \, dy.$$

Consider again two interesting cases. For $k = 0$, we have an equation showing a quantity of all primes up to $x$. In other words, primes counting function can be expressed through the following identity

$$\pi(x) = x H_p(x) - \int_2^x H_p(y) \, dy. \qquad (17)$$

For $k = 1$, we have an equation showing the sum of all primes up to $x$, i.e.:

$$\sum_{p \leq x} p = 2 + 3 + 5 + 7 + \ldots = x^2 H_p(x) - 2 \int_2^x y H_p(y) \, dy.$$

The identity for $H_p(x)$, expressed through set of primes, can be readily found using (13) and substituting $k = 1$:

$$H_p(x) = \frac{1}{x^2} \sum_{p \leq x} p + 2 \int_2^x \frac{1}{y^3} \sum_{p \leq x} p \, dy. \qquad (18)$$

By analogy to (15), RHS of (18) does not include the set of reciprocals. As a result, the algorithm, built on the basis of identity (18), essentially decreases the computation time, especially for the large arguments $x$.

The computation of $H_p(x)$ can also be performed though $\pi(x)$. Rewrite (17) in the following form

$$\pi(x) = \int_{2-\varepsilon}^x y H_p(y)' \, dy. \qquad (19)$$

Rearranging (19) and applying integration by part, we get

$$H_p(x) = \frac{\pi(x)}{x} + \int_2^x \frac{\pi(y)}{y^2} dy.$$

## V. Proof of π-function identity

Dirac delta function approach can also be useful for any kinds of broken-line functions. As an example, we represent a simple proof for the identity relating prime counting function and *Li*-function [11]

$$\pi(x) = Li(x) + \int_2^x \frac{yR(y)'}{\log y} dy - Li(2) + 1, \qquad (20)$$

where $R(x)$ is defined as

$$R(x) = \sum_{p \leq x} \frac{\log p}{p} - \log x. \qquad (21)$$

Using the sampling property (1) of Dirac delta function, $R(x)$ can be expressed as

$$\sum_{p \leq x} \frac{\log p}{p} - \log x = \int_{2-\varepsilon}^x \frac{\log y}{y} \sum_{p \leq x} \delta(y - p) dy - \log x. \qquad (22)$$

The derivative of (22) yields

$$R(x)' = \frac{\log x}{x} \sum_{p \leq x} \delta(x - p) - \frac{1}{x},$$

or

$$\frac{1}{\log x} + \frac{xR(x)'}{\log x} = \sum_{p \leq x} \delta(x - p). \qquad (23)$$

Integrating of (23) results to

$$Li(x) + \int_2^x \frac{yR(y)'}{\log y} dy + const. = \sum_{p \leq x} \chi(x - p).$$

The sum of Heaviside step functions $\sum_{p \leq x} \chi(x - p)$ is actually π-function counting the number of all primes up to $x$. A constant in LHS can be determined from the initial condition at $x = 2$ when π-function is unity. This provides exact formula for prime counting function and completes the proof of (20).

The proof can be obtained even simpler. Consider Mertens' formula providing a relation for sum consisting of the reciprocal primes [12]

$$H_p(x) = \log\log(x) + 1 - \log\log 2 + \int_2^x \frac{R(y)}{y(\log y)^2} dy + \frac{R(x)}{\log x}. \qquad (24)$$

Taking derivative of (24) yields

$$H_p(x)' = \frac{1}{x \log x} + \frac{R(x)'}{\log x}. \qquad (25)$$

Ultimately, substituting (25) into (19) leads to (20). It is proven again.

## VI. Conclusion

We applied Dirac delta function approach to derive identities for the sums of the numbers and their reciprocals. $\delta$-Function is also useful to obtain a simple proof of the identity relating prime counting function and *Li*-function.

## Acknowledgements

Authors express gratitude to Prof. T. D. Radjabov and Dr. E. N. Tsoy for their constructive remarks and suggestions.

## References

[1] P. A. M. Dirac, Proc. R. Soc., Series A, 113 (1927) 621

[2] P. A. M. Dirac, *Principles of Quantum Mechanics*, Oxford University Press, 1930

[3] R. P. Kanwal, *Generalized Functions: Theory and Applications*, Birkhäuser, 2004

[4] P. Antosik, J. Mikusinski, R. Sikorski, *Theory of Distributions. The Sequential Approach,* Elsevier Scientific Publishing Company, 1973

[5] A. I. Khuri , *Advanced Calculus with Applications in Statistics*, Wiley, 2003

[6] A. I. Khuri, Int. J. Math. Educ. Sci. Technol., 35 (2004) 185

[7] C. Guilpin, J. Gacougnolle and Y. Simon, Appl. Num. Math., 48 (2004) 27

[8] M. Pharr and G. Humphreys, *Physically Based Rendering: From Theory to Implementation*, The Morgan Kaufmann Series in Interactive 3D Technology, pp. 279-367, 2004

[9] G. H. Hardy, J. E. Littlewood, and G. Pólya, *Inequalities*, Cambridge University Press, pp. 152-155, 1988

[10] We consider an alternative form of Heaviside function when $\chi(0) \equiv 1$. See, e.g.:

http://en.wikipedia.org/wiki/Heaviside_step_function

[11] R. M. Abrarov, S. M. Abrarov (to be published)

[12] G. Tenenbaum and M. M. France, *The Prime Numbers and Their Distribution*, AMS, p. 22, 2001

---

[1] York University, Toronto, Canada — abrarov@yorku.ca
Dongguk University, Seoul, South Korea — abrarov@dongguk.edu

[2] University of Toronto, Canada — rabrarov@physics.utoronto.ca